\documentclass[12pt]{amsart}

\usepackage{geometry}                
\usepackage{lineno}

\usepackage{mathrsfs}
\usepackage{graphicx}
\usepackage{enumerate}
\usepackage{epstopdf}
\usepackage{amsthm}

\usepackage{dsfont}
\usepackage{amssymb,latexsym,amsmath}
\usepackage{amsfonts}
\usepackage[greek,english]{babel}
\usepackage[iso-8859-7]{inputenc}
\usepackage{bm}

\newtheorem{theorem}{Theorem}
\newtheorem{definition}[theorem]{Definition}
\newtheorem{example}{Example}
\newtheorem{lemma}{Lemma}

\begin{document}
\bibliographystyle{plain}

\title {Phase transition in one-dimensional subshifts}
\author{Nicolai Haydn}

\date{}
\maketitle

\begin{abstract}
\noindent In this note we give simple examples of a one-dimensional mixing subshift 
with positive
topological entropy which have two distinct measures of maximal entropy.
We also give examples of subshifts which have two mutually singular
equilibrium states for H\"{o}lder continuous functions. We also indicate 
how the construction can be extended to yield examples with any number
of measures of maximal entropy of equilibrium states.
\end{abstract}

\section{Introduction}

It is a well known fact that a wide class of uniformly hyperbolic 
transformation that satisfy some irreducibility condition have unique 
measures of maximal entropy. Irreducible subshifts have unique measures 
of maximal entropy if they are for instance
of finite type or sofic (factors of shifts of finite type). More
general hyperbolic systems which satisfy specification have unique
measures of maximal entropy \cite{Bow}. In either case the systems have
a `uniform' mixing property. 
Hofbauer had shown in \cite{Hof} that even subshifts of finite 
type don't necessarily have unique equilibrium states. However
the potentials in the provided examples lack of course regularity 
and one 
of the equilibrium states is supported on a single fixed point. 
On infinite alphabets there are examples of Markov shifts  that
have many Gibbs states although all are non-invariant~\cite{Geo} with only 
one exception.
For a large class of H\"{o}lder continuous functions 
we get here multiple equilibrium states which are supported on
embedded subshifts.

Here we produce a simple example which 
shows that even in the case of an irreducible subshift with positive 
entropy one cannot generally expect to have a unique measure of maximal 
entropy. 
Also Krieger~\cite{Kr} gave an example with a 
non-unique measure of maximal entropy---although not as simple and 
straightforward as this one.
For two-dimensional subshifts of finite type, Burton and Steiff~\cite{B-S, BS2}
have given examples whose measures of maximal entropy form
simplices that consist of more than single points. 

In section~2 we give an example of a subshift which has two 
measures of maximal entropy. We then extend it to examples
which can have any number of measures of maximal entropy
including one with infinitely many. 
In the third section we construct an example that has two 
equilibrium states for H\"older continuous potentials
which we then extend to examples with any number of 
equilibrium states and also one with infinitely many.

\section{Measures of maximal entropy}   

We will first describe an example with two distinct ergodic measures
of maximal entropy and then show how to construct an example with
any number of distinct ergodic measures including one with infinitely 
many.

Let $\nu\ge2$ be some positive integer and let us construct a subshift over 
the symbolset $\mathcal{A}=\{0,1,2,\dots,2\nu\}$. We
shall refer to the symbols $\{1,2,\dots,\nu\}$ as the green symbols 
and denote them by $\mathcal{G}$. Similarly we call $\{\nu+1,\dots,2\nu\}$ 
the yellow symbols and write 
$\mathcal{Y}$ for them. Hence $\mathcal{A} =\{0\}\cup \mathcal{G} \cup \mathcal{Y}$.
Let us denote by $\Sigma_\mathcal{G} =\mathcal{G}^{\bf Z}$ the (green) full $\nu$-shift
over the alphabet $\mathcal{G}$. Similarly we write $\Sigma_\mathcal{Y}$
for the (yellow) full $\nu$-shift $\mathcal{Y}^{\bf Z}$.

Let $\tau>0$. 
The subshift $\Sigma \subset \mathcal{A}^{\bf Z}$ will consist of the union
of the two monochromatic shift spaces  $\Sigma_\mathcal{G}$ and 
$\Sigma_\mathcal{Y}$, to which we add bicoloured sequences in the 
following way: A sequence in $\Sigma$ in which a word $\alpha$ of 
one colour (green or yellow) is followed 
by a word $\beta$ of the other colour (yellow or green) will have a string
of zeroes $\gamma$ separating the two words $\alpha$ and $\beta$. 
The length $|\gamma|$ of that string of zeros is then at least, 
$\tau(a + b)$, where $a = |\alpha|, b = |\beta|$ are the lengths of 
the coloured words to the `left' and `right' of $\gamma$. 
If in particular the word $\alpha$ or $\beta$ is infinitely long then 
it can only be followed by an infinite string of zeros and cannot
`transition' to a string of another colour.

More precisely:

\begin{definition} The shift space $\Sigma \subset \mathcal{A}^{\bf Z}$ 
is defined as follows: We have $x \in \Sigma$ if

\noindent (i) either  $x \in \Sigma_\mathcal{G} \cup \Sigma_\mathcal{Y}$, in which
case we call $x$ a {\em monochromatic} point.

\noindent (ii) or if bi-coloured blocks are of the form:
$$\dots 0g_1g_2\dots g_{a-1}g_a 0^{\lambda}y_1y_2\dots y_{b-s}y_b0\dots$$
or
 $$\dots 0y_1y_2\dots y_{a-1}y_a 0^{\lambda}g_1g_2\dots g_{b-s}g_b0\dots$$
where $\lambda \geq \tau(a + b)$ and $g_i \in \mathcal{G}, y_i \in\mathcal{A}$
($0^{\lambda}$ is a string of zeros of length $\lambda$), $\tau>0$.

\end{definition}

\noindent 
If $\omega = \omega_m \omega_{m+1} \dots \omega_{m+n-1}$ is an allowed word 
of some length $n$ in $\Sigma$, then we denote by $U(\omega)$ the
{\em cylinder set} $\{x \in \Sigma: x_i = \omega_i, m \leq i < m + n \}$.
The shift transformation $\sigma$ on $\Sigma$ is defined in the usual
way by $\sigma(x)_i = s_{i+1}, i \in {\bf Z}$, for $x \in \Sigma$.
Notice that the transition rules are shift invariant. Thus 
$\Sigma$ is indeed a subshift which we endow with the usual topology
generated by the cylinder sets $U(\omega)$, where $\omega$ 
runs over all finite words in $\Sigma$.


If $\tau=0$ then $\Sigma$ is a sofic system in which every string of
zeroes connects to coloured symbols of different colours,
and whose topological entropy is $2\log\nu$ (see lemma 2).

\begin{lemma} The shift transformation $\sigma$ on $\Sigma$ is 
topologically mixing for every $\tau>0$.
\end{lemma}

\noindent {\bf Proof.}
 We have to show that for any two finite words $\omega$
and $\eta$ there exists a number $N$ so that 
$U(\omega) \cap \sigma^n(U(\eta))$ is non-empty for all $n \geq N$.
This is evidently true if $\omega$ and $\eta$ are both monochromatic
of the same colour. If the last symbol of $\eta$ and the 
first symbol of $\omega$ are of different colour then consider the word
$\pi = \eta 0^{\lambda} \omega$. 
Clearly, $\pi$ is an admissable word for every 
$\lambda > \tau(|\eta| + |\omega|)$, 
and the cylinder set 
$U(\pi) \subset U(\omega) \cap \sigma^n(U(\eta))$ is non-empty for all
$n \geq N$, where  $N$ can then be chosen to be equal to
$(1+\tau)(|\eta| + |\omega|)$.

If the last symbol of $\eta$ and the first symbol of $\omega$ are of
the same colour, then let us consider the word 
$\pi = \eta 0^{\kappa} \varepsilon 0^{\lambda} \omega$,
where $\varepsilon$ is a symbol (or some other short word) of
the other colour. Here we have $\kappa \geq \tau(|\eta| + |\varepsilon|)$
and $\lambda \geq \tau(|\omega| + |\varepsilon|)$. Clearly
$U(\pi)$ is non-empty for those choices of $\lambda$ and $\kappa$.
Thus $U(\omega) \cap \sigma^n(U(\eta)) \not= \emptyset$
if $n \geq N = (\tau+1)(|\eta| + 2|\varepsilon| + |\omega|)$.
\hfill$\Box$

\begin{lemma} If $\tau\ge\frac{\log3}{\log\nu}$ then the topological entropy $\,h$ of the shift 
$\sigma: \Sigma \rightarrow \Sigma$ is equal to $\log \nu$.
\end{lemma}

\noindent {\bf Proof.} 
Let us first get a lower bound on the topological entropy $h$ 
of $\Sigma$. For that purpose let us observe that since 
$\Sigma$ contains two full $\nu$-shifts, $\Sigma_\mathcal{Y}$ and $\Sigma_\mathcal{G}$,
the topological entropy of $\Sigma$  must be at least  $\log\nu$. 

Let us now estimate the number of the remaining words of length 
$n$ from above. We have two cases to consider:

\noindent (i) Monochromatic words that might also contain zeroes.
According to our rules such words begin or end with strings
of zeroes. Thus we obtain purely yellow and green words of lengths
$k = 0,\dots,n$ which on at least one side are framed by strings of zeroes. 
Their number turns out to be
$2 \sum_{k=0}^n (n-k)\nu^k\le2n\nu^n\sum_{k=0}^n \nu^{-k}$
which has exponential growth rate $\log\nu$.

\noindent (ii) 
To estimate the number of words of length  $n$  that genuinely contain
symbols of both colours, let us observe that since any such word has at 
least one transition from yellow to green or vice versa, it therefore 
must also contain at least  $n\frac{\tau}{\tau+1}$ zeroes,
that is at most $n'=[\frac{n}{\tau+1}]$ ($[\;]$ denotes integer part)
 coloured symbols (i.e.\ symbols in  $\mathcal{A}\setminus\{0\}=\{1,2,\dots,2\nu\}$). 
The coloured symbols come in monochromatic blocks of 
alternating colour. Denote by $P_{k,\ell}$ the number of possibilities
in which one can arrange $\ell$ symbols in $k$ blocks (separated by
the appropriate number of zeroes), where $k=1,2,\dots,\ell$.
One finds that 
$$
P_{k,\ell}=\left(\begin{array}{c}\ell-1\\k-1\end{array}\right),
$$
which is the number of possibilities of picking the first element
of every block but the very first one.
The number of $n$-words
in $\Sigma$ which contain $\ell$ coloured symbols arranged 
in $k\leq\ell$ monochromatic blocks of alternating colour
 is $2P_{k,\ell}\nu^{\ell}$.

Distributing $\ell$ coloured symbols out of $n' $ in $k$ blocks can
be done in $P_{k,\ell}$ many ways, where $1\leq k\leq\ell\leq n'$. 
This leaves $m=n-(\tau+1)\ell$ zeroes to be distributed on $k+1$ intervals, 
namely the $k-1$ gaps between the blocks of coloured symbols plus
the two ends of the entire word. There are 
$$
Q_{k,m}=\left(\begin{array}{c}m+k\\k\end{array}\right)
$$
many possibilities. By Stirling's formula we have
\begin{eqnarray*}
Q_{k,m}&=&\left(\begin{array}{c}n-(\tau+1)\ell+k\\k\end{array}\right)\\
&\le&\left(\begin{array}{c}n-(\tau+1)\ell+k\\ \frac12(n-(\tau+1)\ell+k)\end{array}\right)\\
&\le& c_12^{n-(\tau+1)\ell+k}\sqrt{n-(\tau+1)\ell+k}\\
&\le&c_12^{n-(\tau+1)\ell+k}\sqrt{n}.
\end{eqnarray*}
We can thus estimate the total number of bi-coloured strings of length $n$ by
\begin{eqnarray*}
Q(n)&\le&2\sum_{\ell=2}^{n'}\nu^\ell \sum_{k=2}^\ell
\left(\begin{array}{c}\ell-1\\k-1\end{array}\right)\left(\begin{array}{c}n-(\tau+1)\ell+k\\k\end{array}\right)\\
&\le&2 c_1\sqrt{n}\sum_{\ell=2}^{n'}\nu^\ell2^{n-(\tau+1)\ell}\sum_{k=2}^\ell
\left(\begin{array}{c}\ell-1\\k-1\end{array}\right)2^k\\
&\le&2 c_1\sqrt{n}\sum_{\ell=2}^{n'}\nu^\ell2^{n-(\tau+1)\ell}3^\ell\\
&\le&c_2\sqrt{n}\left\{\begin{array}{lll}2^n&\mbox{\rm if}&3\nu2^{-1-\tau}<1\\
2^n(3\nu2^{-1-\tau})^\frac{n}{\tau+1}&\mbox{\rm if}&3\nu2^{-1-\tau}\ge1\end{array}\right.
\end{eqnarray*}
Thus 
$$
\lim_{n\rightarrow\infty}\frac1n\log Q(n)\le\max \left(\log2,\frac{\log3\nu}{\tau+1}\right)
$$
which is $\le\log\nu$ if $\tau\ge\frac{\log3}{\log\nu}$.

 \hfill$\Box$

\begin{example}\label{singular.measures}
 If $\tau\geq\frac{\log3}{\log\nu}$, then there are two mutually 
singular measures of maximal entropy on $\Sigma$.
\end{example}

\noindent {\bf Proof.} 
It is well known \cite{Wal} that the measure of maximal entropy on 
the full 
$\nu$-shift $\{1,\dots,\nu\}^{\bf Z}$ is the Bernoulli measure with the 
probability vector $(\frac1{\nu},\dots,\frac1{\nu})$. It's metric entropy 
is $\log{\nu}$. Let us define on $\Sigma$ a Bernoulli measure
$\mu_\mathcal{Y}$ which is given by the probability vector 
$(0,\frac1{\nu},\dots,\frac1{\nu},0,\dots,0)$ and the measure 
$\mu_\mathcal{G}$ 
with the probabilities $(0,0,\dots,0,\frac1{\nu},\dots,\frac1{\nu})$.
Evidently both measures are shift invariant and have metric
entropies $\log{\nu}$ which by Lemma~2 is the topological entropy
of $\Sigma$. Hence $\mu_\mathcal{Y}$ and $\mu_\mathcal{G}$ 
(as well as all their linear combinations) are
distinct ergodic measures of maximal entropy for the subshift $\Sigma$.
\hfill$\Box$

\subsection{Example with $L\ge2$ distinct ergodic measures of maximal entropy}
Let $\nu\ge2$ and integer and consider the alphabet 
$\mathcal{A}=\{0,1,2,\dots,L\nu\}$. We label the embedded full $\nu$-shift spaces
by $\Sigma_j=\{(j-1)\nu+1,(j-1)\nu+2,\dots,j\nu\}^\mathbb{Z}$, $j=1,2,\dots, L$.
A transition from a word $\alpha$ in $\Sigma_j$ happens only to a word $\beta$ 
of a `different colour' if it is separated by a string of zeros of length $\ge\tau(|\alpha|+|\beta|)$
and $\beta$ lies either in $\Sigma_{j-1}$ or $\Sigma_{j+1}$. 
One sees, as was proven in Lemma~1 for $L=2$, that the resulting shift space 
$\Sigma\subset\{0,1,2,\dots,L\nu\}^\mathbb{Z}$ is topologically mixing. 

\begin{example}
If $\nu\ge2$ and $\tau\ge\frac{\log5}{\log\nu}$ then there are $L\ge2$ distinct ergodic measures
of maximal entropy on $\Sigma$.
\end{example}

\noindent {\bf Proof.} The topological entropy of $\Sigma$ is at least $\log \nu$ as this is the topological
entropy of the embedded subshifts $\Sigma_j$. To get an upper bound on the 
topological entropy of $\Sigma $ we estimate the number $Q(n)$ of words of length $n$
that are not `monochromatic', i.e.\ don't belong to a single $\Sigma_j$. We proceed as 
before and obtain:
$$
Q(n)\le L\sum_{\ell=2}^{n'}\nu^\ell\sum_{k=2}^\ell
\left(\begin{array}{c}\ell-1\\k-1\end{array}\right)2^{k-1}\left(\begin{array}{c}m+k\\k\end{array}\right),
$$
where $m=n-[(\tau+1)\ell]$,  $\ell$ counts the number of `coloured' symbols 
(i.e.\ symbols in $\mathcal{A}\setminus\{0\}$), $k-1$ is the number of changes of colour 
that is
where there is a transition from a word in some $\Sigma_j$ to $\Sigma_{j\pm1}$
(the factor $2^{k-1}$ accounts for number of possible switches in colour). 
As before we thus get (with $n'=[n/(\tau+1)]$)
\begin{eqnarray*}
Q(n)&\le&L\sqrt{n}c_1\sum_{\ell=2}^{n'}\nu^\ell2^{n-(\tau+1)\ell} \sum_{k=2}^\ell
\left(\begin{array}{c}\ell-1\\k-1\end{array}\right)2^k2^{k-1}\\
&\le&2L \sqrt{n}c_1\sum_{\ell=2}^{n'}2^n(5\nu2^{-\tau-1})^\ell\\
&\le&c_2\sqrt{n}\left\{\begin{array}{lll}2^n&\mbox{\rm if}&5\nu2^{-\tau-1}<1\\
2^n(5\nu2^{-\tau-1})^{\frac{n}{\tau+1}}&\mbox{\rm if}&5\nu2^{-\tau-1}\ge1\end{array}\right..
\end{eqnarray*}
Consequently we get that 
$$
Q(n)\le c_2\sqrt{n}\,2^n\left(1+(5\nu2^{-\tau-1})^{\frac{n}{\tau+1}}\right)
$$
and thus 
$$
\lim_{n\rightarrow\infty}\frac1n\log Q(n)\le\max\left(\log2,\frac{\log5\nu}{1+\tau}\right).
$$
Hence, if $\tau\ge\frac{\log5}{\log\nu}$ then $h_{\rm top}(\Sigma)=\log \nu$.
It is now clear how we get the $L$ distinct ergodic measures of maximal entropy:
they are the $L$ Bernoulli measures $\mu_j$  with weights 
$p^j_k=\frac1\nu$ for $k=(j-1)\nu+1,(j-1)\nu+2,\dots,j\nu$ and $p^j_k=0$ for all other $k$
($j=1,2,\dots, L$). These measures all mutually singular and have metric entropy $\log \nu$ i.e.\ they are
measures of maximal entropy.
\hfill$\Box$

\subsection{Example with infinitely many distinct ergodic measures of maximal entropy}
Let $\nu\ge2$.
We now let $L$ go to infinity and obtain a subshift $\Sigma$ over the alphabet
$\mathcal{A}=\mathbb{N}_0=\{0,1,2,\dots\}$.  There are infinitly many embedded full $\nu$-shift
$$
\Sigma_j=\{(j-1)\nu+1,(j-1)\nu+2,\dots,j\nu\}^\mathbb{Z},
$$ 
for $j=1,2,\dots$.
As before we only transition from a word $\alpha$ in $\Sigma_j$ a word $\beta$ 
in $\Sigma_{j\pm1}$ by keeping them separted by a string of zeros of length $\ge\tau(|\alpha|+|\beta|)$.
Again the resulting shift space $\Sigma$ is topologically mixing but it is not a compact space. 
However we can consider its one-point-compactification $\bar\Sigma=\Sigma\cup\{\infty\}$
and a suitable metric which in particular implies that neighbourhoods of $\infty$ contain the 
embedded subshifts $\Sigma_j$ for large enough $j$. The topological entropy can then 
be defined by $(\varepsilon,n)$-separated sets $E_{\varepsilon,n}$ (see e.g.~\cite{Wal}) and 
is given by 
$h_{\rm top}(\Sigma)=\lim_{\varepsilon\rightarrow0^+}\lim_{n\rightarrow\infty}\frac{\log|E_{\varepsilon,n}|}n$.

\begin{example}
If $\nu\ge2$ and $\tau\ge\frac{\log5}{\log\nu}$ then there are infinitely many distinct ergodic measures
of maximal entropy on $\Sigma$.
\end{example}

\noindent {\bf Proof.} Evidently $h_{\rm top}(\Sigma)\ge\log\nu$. Moreover if we denote by
$R_L(n)$ the number of words of length $n$ in $\Sigma$ whose first coordinate lie in
$\mathcal{A}_L=\{0,1,\dots,L\nu\}$ then we get that 
$h_{\rm top}(\Sigma)=\lim_{L\rightarrow\infty}\lim_{n\rightarrow\infty}\frac{\log R_L(n)}n$.
We have that 
$$
R_L(n)\le L\nu^n+Lc_3\sqrt{n}\left\{\begin{array}{lll}2^n&\mbox{\rm if}&5\nu2^{-\tau-1}<1\\
2^n(5\nu2^{-\tau-1})^{\frac{n}{\tau+1}}&\mbox{\rm if}&5\nu2^{-\tau-1}\ge1\end{array}\right..
$$
and consequently 
$$
\lim_{n\rightarrow\infty}\frac1n\log R_L(n)\le\max\left(\log \nu,\frac{\log5\nu}{1+\tau}\right).
$$
Now, letting $L$ go to infinity we obtain 
$h_{\rm top}(\Sigma)\le\max\left(\log \nu,\frac{\log5\nu}{1+\tau}\right)$
which is equal to $\log\nu$ if $\tau\ge\frac{\log5}{\log\nu}$ then $h_{\rm top}(\Sigma)=\log \nu$.
The equal weight Bernoulli measures on the spaces $\Sigma_j$ are all measures
of maximal entropy and are mutually singular.
\hfill$\Box$


\section{Equilibrium states}

\noindent
Let us construct now subshifts with  non unique equilibrium states for H\"older continuous
potentials.
We say a function $f:\Sigma\rightarrow\mathbb{R}$ is H\"older continuous on the
shift space $\Sigma$ if there exists a $\vartheta\in(0,1)$ and a constant $C$ so that
for every $n>0$ one has
$$
\sup|f(\vec{x})-f(\vec{y})|\leq C\vartheta^n
$$
where the supremum is over all pairs of sequences 
$\vec{x}=\cdots x_{-1}x_0x_1x_2\dots,
\vec{y}=\cdots y_{-1}y_0x_1y_2\dots\in\Sigma$
for which $x_i=y_i\forall |i|<n$ and either
$x_{-n}\not=y_{-n}$ or $x_n\not=y_n$.

\subsection{Example with two distinct equilibrium states}
Let $\tau>0$ and $\nu\ge3$ and integer. Then we put as above
$\mathcal{A}=\{0\}\cup\mathcal{G}\cup\mathcal{Y}=\{0,1,\dots,2\nu\}$ and
define as above the subshift $\Sigma\subset\mathcal{A}^\mathbb{Z}$
which has the two embedded full $\nu$-shifts $\Sigma_*$,  $*=\mathcal{G},\mathcal{Y}$
and where the transitions are strings of zeros of lengths $\ge\tau(|\alpha|+|\beta|)$
with $\alpha,\beta$ being the adjacent monochromatic words of different colours.

Let $f$ be a real valued H\"older continuous function on the full 
two-shift $\Sigma_\mathcal{G}$
 so that $P_\mathcal{G}(f)>\sup f+\log 2$, where $P_\mathcal{G}(f)$ is the pressure 
of $f$ on the full shift $\Sigma_\mathcal{G}$ (see \cite{Wal}). For $\nu\ge3$ this is
possible since $P_\mathcal{G}(f)-\sup f$ is bounded above by $\log \nu$ and 
equal to $\log \nu$ for constant $f$ and the pressure function is continuous.

Let us define a H\"older continuous function $g$ on $\Sigma$ in two 
steps. We first extend the function to $\Sigma_\mathcal{Y}$ by identification:
that is $g(\vec{y})=f(\vec{x})$ for $\vec{x}\in\Sigma_\mathcal{G}, \vec{y}\in\Sigma_\mathcal{Y}$
for which $y_i=x_i+\nu\forall i$. Now we extend $g$ to bi-coloured 
words in the following way: If $x_0=0$ then we put $g(\vec{x})=\sup f$.
If $x_0\not=0$ and $\vec{x}\in\Sigma$ is a sequence which contains
symbols of both colours then let $x_{-n}\cdots x_n$ be the 
longest monochromatic
string of either green or yellow symbols (that is $x_{-n}\cdots x_n$ is 
a word in either $\Sigma_\mathcal{G}$ or $\Sigma_\mathcal{Y}$ and either $x_{-n-1}=0$
or $x_{n+1}=0$). We then pick a monochromatic point 
$\vec{y}\in U(x_{-n}\cdots x_n)\subset\Sigma_*$
($*=\mathcal{G},\mathcal{Y}$) and define $g(\vec{x})=f(\vec{y})$.
In this way $g$ is a H\"older continuous function on $\Sigma$.
Let us now prove that the pressure of $g$ (on $\Sigma$) is  in fact
equal to $P_\mathcal{G}(f)$ (which equals $P_\mathcal{Y}(g)$ for a suitable 
choice of $\tau$.

\begin{lemma}
If $\tau\geq \frac{\log4\nu}{P_\mathcal{G}(f)-\sup f} -1$, then the pressure
$P(g)$ of $g$ on $\Sigma$ is equal to $P_\mathcal{G}(f)$.
\end{lemma}

\noindent {\bf Proof.} 
Denote by $W_n$ the collection of $\Sigma$-words
$\alpha=\alpha_1\cdots \alpha_n$ of lengths $n$.
The pressure of $g$ is then given by the exponential 
growth rate of the partition function
$$
Z_n=\sum_{\alpha\in W_n} e^{g^n(x(\alpha))},
$$
where $x(\alpha)$ is an arbitrary point in the cylinder set 
$U(\alpha)\subset\Sigma$, and 
$g^n=g+g\sigma+g\sigma^2+\cdots+g\sigma^{n-1}$ is the $n$-th
ergodic sum of $g$.

If in the above sum we restrict to monochromatic words $\alpha$ of
either colour, then we obtain the green or yellow
partition function $Z^*_n=\sum_{\alpha\in W^*_n} e^{g^n(x(\alpha))}$ ($*=\mathcal{G},\mathcal{Y}$)
whose growth rates are exactly $P_\mathcal{G}(f)$. We denoted by
$W^*_n$ the set of monochromatic words $\alpha=\alpha_1\alpha_2\cdots\alpha_n$
in $\Sigma_*$.

Let us note that the growth rate of 
$\sum_{k=1}^nZ^*_k$ is also $P_\mathcal{G}(f)$.  Hence, if we denote by $W'_n$
those words in $W_n$ that contain symbols of exactly one colour and
otherwise at least one zero, then it follows that growth rate of 
$Z_n=\sum_{\alpha\in W'_n} e^{g^n(x(\alpha))}$
is also given by $P_\mathcal{G}(f)$.

It thus remains to show that the exponential growth rate of 
$\sum_{\alpha\in W''_n} e^{g^n(x(\alpha))}$
is $\leq P_\mathcal{G}(f)$, where $W''_n$ denotes the genuinely bi-coloured 
words in $W_n$ (i.e.\ those words that contain green and yellow symbols).
In Lemma~2 we estimated 
$$
|W''_n|=Q(n)\leq c_2\sqrt{n}2^n\left(1+ (3\nu2^{-1-\tau})^{n/(\tau+1)}\right),
$$
and since $\sup g=\sup f$, we obtain
$$
Z''_n=\sum_{\alpha\in W''_n} e^{g^n(x(\alpha))}
\leq |W''_n|e^{n\sup g}
\leq c_2e^{n\sup f}\sqrt{n}\,2^n\left(1+ (3\nu2^{-1-\tau})^{n/(\tau+1)}\right).
$$
Thus 
$$
\lim_{n\rightarrow\infty}\frac{\log Z_n''}n\le\sup f+\max\left(\log2,\frac{\log3\nu}{\tau+1}\right)
$$
which is $\leq P_\mathcal{G}(f)$ if $\tau+1\geq \frac{\log3\nu}{P_\mathcal{G}(f)-\sup f}$.
\hfill$\Box$

\begin{example} Let $\nu\ge3$ and $\tau\geq \frac{\log3\nu}{P_\mathcal{G}(f)-\sup f} -1$.
Then there are two mutually 
singular equilibrium states for $g$ on $\Sigma$.
\end{example}

\noindent {\bf Proof.} For $\nu\ge3$ there are many potentials 
$f:\Sigma_\mathcal{G}\rightarrow\mathbb{R}$
for which $P_\mathcal{G}(f)\ge\sup f+\log2$. Then for the extension 
$g:\Sigma\rightarrow\mathbb{R}$
there is an equilibrium state on each of the embedded full $\nu$-shifts
 $\Sigma_\mathcal{G}$ and $\Sigma_\mathcal{Y}$ both of which,
by Lemma~4, 
assume the pressure $P(g)=P_\mathcal{G}(f)$ in the variational principle. 
These two measures are both invariant and mutually singular. 
\hfill$\Box$

\subsection{Example with $L\ge2$ distinct equilibrium states}
For $\nu\ge3$, an integer, and construct as above a subshift $\Sigma$
over the alphabet $\mathcal{A}=\{0,1,2,\dots,L\nu\}$ so that it contains the embedded 
full $\nu$-shift spaces $\Sigma_j=\{(j-1)\nu+1,\dots,j\nu\}^\mathbb{Z}$, $j=1,2,\dots, L$.
Again we allows for transitions from words $\alpha$ in some $\Sigma_j$  to a word $\beta$ 
in $\Sigma_{j\pm1}$ that have to separated by a string of zeros of length $\ge\tau(|\alpha|+|\beta|)$.

Let $f$ be a real valued H\"older continuous function on the full 
two-shift $\Sigma_1$
 so that $P_1(f)>\sup f+\log 2$, where $P_1(f)$ is the pressure 
of $f$ on the full shift $\Sigma_1$. 

Let us define a H\"older continuous function $g$ on $\Sigma$ in two 
steps. We first extend the function from $\Sigma_1$ to $\Sigma_j$ for 
$j=2,\dots,L$ by identification: If $\vec{x}\in\Sigma_1, \vec{y}\in\Sigma_j$
so that $y_i=x_i+(j-1)\nu\forall i$, then we put $g(\vec{y})=f(\vec{x})$.
Then one  extends $g$ to bi-coloured 
words as before, namely we put $g(\vec{x})=\sup f$ if $x_0=0$ and if 
$x_0\not=0$ for a non-monochromatic $\vec{x}\in\Sigma$ then we put 
$g(\vec{x})=f(\vec{y})$ where $\vec{y}$ is a monochromatic point that
realises the longest symmetric monochromatic word of $\vec{x}$. 

\begin{example}
If $\tau\geq \frac{\log4\nu}{P_\mathcal{G}(f)-\sup f} -1$, then the pressure
$P(g)$ of $g$ on $\Sigma$ is equal to $P_\mathcal{G}(f)$.
\end{example}

\noindent {\bf Proof.} 
We have to find the exponential growth rate $P(g)$ of the partition function 
$
Z_n=\sum_{\alpha\in W_n} e^{g^n(x(\alpha))},
$
where  $W_n$ are the $\Sigma$-words of lengths $n$. Clearly 
$P(g)\ge P_1(f)$.

In order to estimate the exponential growth rate of 
$Z''_n=\sum_{\alpha\in W''_n} e^{g^n(x(\alpha))}$ where $W''_n$ are the
 genuinely multi-coloured words in $W_n$ we get as above that
$$
|W''_n|=Q(n)\leq c_2\sqrt{n}2^n\left(1+ (5\nu2^{-1-\tau})^{n/(\tau+1)}\right),
$$
and therefore
$$
Z''_n=\sum_{\alpha\in W''_n} e^{g^n(x(\alpha))}
\leq |W''_n|e^{n\sup g}
\leq c_2e^{n\sup f}\sqrt{n}2^n\left(1+ (5\nu2^{-1-\tau}))^{n/(\tau+1)}\right)
$$
which yields
$$
\lim_{n\rightarrow\infty}\frac{\log Z_n''}n\le\sup f+\max\left(\log2,\frac{\log5\nu}{\tau+1}\right).
$$
This is  $\leq P_\mathcal{G}(f)$ if $\tau+1\geq \frac{\log5\nu}{P_\mathcal{G}(f)-\sup f}$.
\hfill$\Box$

\begin{lemma} Let $\nu\ge3$ and $\tau\geq \frac{\log5\nu}{P_1(f)-\sup f} -1$.
Then there are $L\ge2$ mutually 
singular equilibrium states for $g$ on $\Sigma$.
\end{lemma}

\noindent {\bf Proof.} For $\nu\ge3$ there are many potentials 
$f:\Sigma_\mathcal{G}\rightarrow\mathbb{R}$
for which $P_\mathcal{G}(f)\ge\sup f+\log2$. Then for the extension 
$g:\Sigma\rightarrow\mathbb{R}$
there is an equilibrium state on each of the embedded full $\nu$-shifts
 $\Sigma_\mathcal{G}$ and $\Sigma_\mathcal{Y}$ both of which,
by Lemma~4, 
assume the pressure $P(g)=P_\mathcal{G}(f)$ in the variational principle. 
These two measures are both invariant and mutually singular. 
\hfill$\Box$

\subsection{Example with infinitely many distinct equilibrium states}
Let $\nu\ge3$ and construct the subshift $\Sigma\subset\mathbb{N}_0^\mathbb{Z}$
as was done above. A H\"older continuous function $f$ on $\Sigma_1$ which 
satisfies $P_1(f)\ge\sup f+\log2$ is then as above extended to a H\"older function
$g$ on $\Sigma$ so that it is a `copy' on each of the embedded full shifts
$\Sigma_j$, $j=1,2,\dots$. Moreover we can extend $g$ to the compact space
$\bar\Sigma$ by defining $g(\infty)=\sup f$. If we denote by $\mu_j$ the equilibrium
state on $\Sigma_j$ for the potential $g|_{\Sigma_j}$ then $P_j(g|_{\Sigma_j})=P_1(f)=P(g)$
for all $j$. We thus obtain:

\begin{example}
If $\nu\ge3$ and $\tau\ge\frac{\log5}{\log\nu}$ and $\tau\geq \frac{\log5\nu}{P_1(f)-\sup f} -1$
then there are infinitely many distinct equilibrium states $\mu_j$ on $\Sigma$ for $g$ which are 
mutually singular.
\end{example}

\end{document}